\documentclass[12pt,french]{article}
\usepackage{amssymb,amsfonts,a4,epsfig,isolatin1,babel}
\newlength{\Taille}

\newcommand{\deriv}[2]{\frac{\partial #1}{\partial #2}}
\newcommand{\ddt}{\frac{d}{dt}}

\newcommand{\pscal}[2]{\langle #1,#2\rangle}
\newcommand{\ham}[1]{\mathcal{X}_{#1}}
\newcommand{\Lie}[1]{\mathfrak{#1}}

\newcommand{\phy}{\varphi}

\newcommand{\om}{\omega}
\newcommand{\al}{\alpha}
\newcommand{\la}{\lambda}

\newcommand{\fleche}{\rightarrow}

\newcommand{\restr}{\upharpoonright}

\newcommand{\demo}{\noindent{\bf D\'emonstration: }}
\newcommand{\rem}[1][]{\noindent {\sl ${\mathcal{R}}$emarque #1 :~}}

\newcommand{\flechemap}[1]{
        \settowidth{\Taille}{\mbox{$ ~ #1 ~$}}
        \stackrel{ ~ #1 ~ }{
        \rule{\the\Taille}{0.1ex}\raisebox{-.5ex}{$\!\fleche$}}}

\newcommand{\cqfd}{\hfill $\square$}

\newcommand{\Cinf}{C^\infty}

\newcommand{\opif}{op\'erateur int\'egral de Fourier}

\newcommand{\opd}{op\'erateur pseudo-diff\'erentiel}
\newcommand{\opds}{op\'erateurs pseudo-diff\'erentiels}

\newcommand{\diffeo}{diff\'eo\-mor\-phis\-me}

\newcommand{\RM}{\mathbb{R}}
\newcommand{\ZM}{\mathbb{Z}}

\newcommand{\NM}{\mathbb{N}}
\newcommand{\CM}{\mathbb{C}}

\newcommand{\T}{\mathbb{T}}

\newcommand{\PM}{\mathbb{P}}

\newcommand{\A}{\mathcal{A}}

\newtheorem{theo}{Th\'eor\`eme}[section]
\newtheorem{defi}{D\'efinition}[section]
\newtheorem{prop}[theo]{Proposition}
\newtheorem{lemm}[theo]{Lemme}
\newtheorem{coro}[theo]{Corollaire}

\newcommand{\He}{\mathcal{H}}

\newcommand{\ad}[1]{\mathrm{ad}_{#1}}

\title{\textsf{Formes normales
semi-classiques des systèmes complètement intégrables au voisinage
d'un point critique de l'application moment}}

\author{V\~u Ng\d oc San\footnote{Institut Fourier UMR5582, B.P. 74, 38402
    Saint-Martin d'H\`eres, France. Actuellement d\'etach\'e \`a~: Mathematics
Institute, Budapestlaan 6, University of Utrecht, 3508 TA Utrecht, the
Netherlands. e-mail:\texttt{vu-ngoc@math.ruu.nl}}} 

\date{}

\begin{document}
\maketitle

\begin{center}
  \textbf{Abstract}\\
  ~\\
  \parbox{10cm}{The semi-classical study of a 1-dimensional Schrödinger
    operator near a non-degenerate maximum of the potential has lead
    Colin de Verdière and Parisse to prove a microlocal normal form
    theorem for any 1-dimensional pseudo-differential operator with
    the same kind of singularity. We present here a generalization of
    this result to pseudo-differential integrable systems of any
    finite degree of freedom with a Morse singularity. Our results are
    based upon Eliasson's study of critical integrable systems.}
\end{center}

\noindent\textbf{Motivation}

L'étude des propriétés spectrales de l'opérateur de Schrödinger
$-\frac{h^2}{2}\Delta + V$ en régime semi-classique ($h\fleche 0)$ en
fonction de la forme du potentiel $V$ pose naturellement le problème
de savoir dans quelle mesure, étant donné un  $h$-\opd\  $P(h)$, son
symbole principal $p$ influe sur la nature des solutions microlocales de
l'équation~:
 \begin{equation}
\label{schrodinger}
 Pu = O(h^\infty). 
\end{equation}
Dans le cas où $P$ est auto-adjoint à variétés caractéristiques
compactes (par exemple, l'opérateur de Schrödinger sur $\RM$ avec un
potentiel tendant vers l'infini en $\pm\infty$) , on a ainsi accès au
comportement microlocal du spectre semi-classique. Ce problème est
discuté par exemple dans \cite{helffer-robert,helffer-sjostrand} et
les références qui y sont citées.  

Au voisinage de \emph{points réguliers} de $p$, la théorie est bien
connue et ne réserve plus guère de surprise. En effet, $P$ se conjugue
par un \opif\ à un opérateur de dérivation
$\frac{h}{\sqrt{-1}}\deriv{}{x}$. Cette forme normale, qui n'est
qu'une autre formulation des solutions WKB, permet de trouver les
conditions, dites ``de Bohr-Sommerfeld'', sous lesquelles l'équation
(\ref{schrodinger}) admet des solutions $L^2$ (voir par exemple
\cite{duistermaat-oscillatory}).

La présence de \emph{points critiques}, et en particulier instables,
rend les choses plus intéressantes, comme en témoignent par exemple
les articles \cite{brummelhuis,sjostrand,colin-p,colin-p2}. Néanmoins,
comme l'ont remarqué les auteurs de \cite{sjostrand,colin-p}, on peut
encore trouver une \emph{forme normale} au voisinage du point
critique, pourvu que celui-ci soit non-dégénéré. Y.Colin de Verdière
et B.Parisse montrent dans \cite{colin-p2} comment en déduire une
version singulière des conditions de Bohr-Sommerlfeld qui permet une
description très précise du spectre semi-classique.

Le but de cet article est de généraliser le théorème de forme normale
microlocale à des systèmes complètement intégrables d'\opds\ en
dimension quelconque. On obtient des analogues semi-classiques étroits
des résultats classiques d'Eliasson. Lorsque la singularité est de
type elliptique, nos résultats sont à rapprocher des formes normales
de Birkhoff semi-classiques de Sjöstrand \cite{sjostrand2} (on obtient
des résultats plus forts, mais les hypothèses le sont aussi).

Dans un autre travail \cite{san2}, nous montreront comment ces
résultats permettent, via des conditions de Bohr-Sommerfeld
singulières, de décrire très précisément le \emph{spectre conjoint} de
ces systèmes complètement intégrables. 

Mentionnons aussi que Bleher, Kosygin et Sinai \cite{bleher} arrivent
à une description spectrale des systèmes de Liouville (qui
représentent une classe générale de systèmes intégrables en dimension
2), par une approche différente.

\section{Introduction}

Sur une variété symplectique $(M,\om)$ de dimension $2n$, un système
comp\-lè\-te\-ment intég\-rable est la donnée de $n$ fonctions
$f_1,\dots,f_n$ en involution pour le crochet de Poisson défini par la
structure symplectique, et dont les différentielles sont presque
partout indépendantes. Si $H$ est un hamiltonien appartenant à
l'algèbre des fonctions $\Cinf$ fonctionnellement engendrées par
$f_1,\dots,f_n$, il est dit complètement intégrable et les $f_i$ sont
des intégrales premières du mouvement. Un tel système définit une
action hamiltonienne locale de $\RM^n$ dans $M$, d'application
moment~:
\[ F~:~ M\ni m \mapsto (f_1(m),\dots,f_n(m))\in \RM^n,\]
et dont les orbites sont génériquement les fibres lagrangiennes
 $F^{-1}(c)$, $c\in\RM^n$. Plus précisément, on sait (théorème
d'Arnold-Liouville \cite{duistermaat}) que si $c$ est une valeur
régulière de $F$, et si  $F^{-1}(c)$ est compacte, alors au voisinage
de $c$, les ensembles de niveau de $F$ sont des tores lagrangiens qui
s'écrivent $\xi=cst$ dans des bonnes coordonnées symplectiques
$(x,\xi)$.
En outre, dans ces coordonnées --  dites ``actions ($\xi$) - angles
($x$)'' -- le mouvement dans chaque tore est linéaire.

Supposons maintenant qu'on se donne $n$ \opds\ $P_1$, $\dots$, $P_n$
sur une variété $X$ de dimension $n$, qui commutent deux-à-deux, et
dont les symboles principaux forment un système complètement
intégrable sur $T^*X$.  En dehors des points critiques de
l'application moment, l'usage des coordonnées actions-angles
classiques permet alors de se ramener au cas où les opérateurs sont
simplement les $\frac{h}{\sqrt{-1}} \deriv{}{x_i}$. Une telle forme
normale permet d'écrire les conditions de Bohr-Sommerfeld sous
lesquelles on peut résoudre simultanément les équations
$P_iu=O(h^\infty)$ et qui sous certaines hypothèses, donnent la forme
du \emph{spectre conjoint} des $P_i$.  Comme le cas unidimensionnel
mentionné plus haut, ceci est connu depuis longtemps, au moins des
physiciens quantiques. Les preuves mathématiques rigoureuses de ce
phénomène, comprenant les asymptotiques com\-plètes des valeurs propres,
sont données par les travaux d'Anne-Marie Charbonnel \cite{charbonnel}
et Yves Colin de Verdière \cite{colinII}.

La question qui nous intéresse est donc de savoir ce que devient le
système au voisinage d'un point critique de $F$. Au niveau classique,
on est amené à se placer dans l'hypothèse d'un point critique ``de
Morse'', ou ``non-dégénéré'', en un sens précisé en section
\ref{classique}.

Le résultat principal (théorèmes \ref{TQ3} et \ref{TQ2}) peut, en
gros, se formuler ainsi~:

\noindent ``l'algèbre engendrée par les opérateurs $P_i$ est, à $O(h^\infty)$
près, conjuguée par un \opif\ à une algèbre standard ne dépendant que
des dérivées secondes des symboles principaux $p_i$ au point
critique.''

\noindent C'est une 
généralisation de \cite{colin-p}, où le résultat est donné en
dimension 1. On verra cependant que, comme le théorème d'Eliasson le
laissait prévoir, le cas multi-dimensionnel introduit une
subtilité, due au fait que l'algèbre standard en question peut avoir
une formulation légèrement différente selon que les fibres
lagrangiennes sont localement connexes ou non. On en donnera malgré
tout une description complète dans le cas général (proposition
\ref{commutant-quantique}).

\section{Formes normales classiques}
\label{classique}
Le but de cette section est de présenter les théorèmes de formes
normales classiques, qui sont essentiellement dus à Eliasson
\cite{eliasson-these}, mais aussi de mettre clairement en valeur les
ingrédients essentiels pour la quantification semi-classique de la
section suivante, comme le \emph{commutant classique} (proposition
\ref{commutant-classique}) et le \emph{lemme de Poincaré critique}.

Les résultats présentés ici sont locaux, au voisinage du point
critique.  Il est à noter néanmoins que l'étude topologique globale de
tels systèmes complètement intégrables à singularités, initiée par
Fomenko, est de mieux en mieux comprise; à ce sujet, se référer aux
dernières mises au point de Nguyen Tiên Zung \cite{zung}.

\subsection{Le théorème d'Eliasson}
Soit $(M,\om)$ une variété symplectique de dimension $2n$, et $m\in
M$.  Le crochet de Poisson munit l'espace des fonctions $\Cinf$ sur
$M$ d'une structure d'algèbre de Lie. On note $\He$ l'homomorphisme de
la sous-algèbre de Lie des fonctions critiques en $m$ sur l'espace
${\mathcal Q}(2n)$ des formes quadratiques sur $\RM^{2n}=\{(x,\xi)\}$,
qui à $f$ associe sa hessienne. La structure d'algèbre de Lie de
${\mathcal Q}(2n)$ est aussi donnée par le crochet de
Poisson. ${\mathcal Q}(2n)$ est ainsi canoniquement isomorphe à
$\Lie{sp}(2n)$.

Soit $(f_1,\dots,f_n)$ un système complètement intégrable sur $M$. On
dit que $m$ est un point critique du système s'il est critique pour
l'application moment $F$. Dans ce travail, on considèrera toujours que
la singularité est de codimension maximale, au sens où chaque $f_i$
est critique en $m$. On supposera aussi que $f_i(m)=0$, ce qui n'ôte
aucune généralité.

À un tel système complètement intégrable singulier, on peut associer
une sous-algèbre réelle $\mathfrak{c}_F$ de ${\mathcal Q}(2n)$, à
savoir la sous-algèbre engendrée par
\[\{ \He(f_1),\dots,\He(f_n)\}.\] 
On remarque que $\mathfrak{c}_F$ est toujours abélienne.
\begin{defi}
  Un système complètement intégrable singulier d'application moment $F$
  est dit \emph{non-dégénéré} au point $m$ si $\mathfrak{c}_F$ est une
  sous-algèbre de Cartan de ${\mathcal Q}(2n)$.
\end{defi}
Rappelons qu'on appelle sous-algèbre de Cartan d'une algèbre
de Lie semi-simple $\Lie{a}$ une sous-algèbre de Lie $\Lie{h}$ commutative
maximale, et qui vérifie en outre la propriété suivante~:
\begin{equation}
\label{semisimple} \forall H\in \Lie{h},~~ \ad{H} 
\textrm{ est un endomorphisme semi-simple de } \Lie{a}.
\end{equation}
Depuis les travaux de Williamson \cite{williamson}, on sait classifier
les sous-algèbres de Cartan réelles de ${\mathcal Q}(2n)$ modulo
conjugaison par un symplectomorphisme. Le résultat est le suivant~:
\begin{theo}[Williamson]
Soit   $\Lie{a}$ une sous-algèbre de Cartan réelle de ${\mathcal
  Q}(2n)$. Il existe des coordonnées symplectiques linéaires
  $(x_1,\dots,x_n,\xi_1,\dots,\xi_n)$ sur $\RM^{2n}$, et une base
  $q_1,\dots,q_n$ de $\Lie{a}$ telle que chaque $q_i$ ait l'une des
  trois formes suivantes~:
\begin{itemize}
        \item  $q_i = x_i\xi_i$ (singularité \emph{hyperbolique})
        \item  $q_i = x_i^2+\xi_i^2$ (singularité \emph{elliptique})
        \item  $q_i = x_i\xi_{i+1} - x_{i+1}\xi_i$, auquel cas on
demande que $q_{i+1} = x_i\xi_i + x_{i+1}\xi_{i+1}$ (singularité
\emph{focus-focus}).
\end{itemize}
\end{theo}
Suivant Eliasson \cite{eliasson-these}, on appellera $(q_1,\dots,q_n)$
une \emph{base standard} de $\Lie{a}$, et suivant Nguyên Tiên Zung, on
dira que $\Lie{a}$ est de \emph{type} $(m_e,m_h,m_f)$, avec
$m_e+m_h+2m_f=n$, si une base standard contient $m_e$ éléments de type
elliptique, $m_h$ éléments hyperboliques, et $m_l$ paires de type
focus-focus.

\rem la condition (\ref{semisimple}) est nécessaire; par exemple, pour
$\Lie{a}=\Lie{sp}(2,\RM)=\Lie{sl}(2,\RM)$, la sous-algèbre engendrée
par $q=\xi^2$ est commutative maximale mais ne relève pas de la
classification ci-dessus. En effet, $\ad{q}=2\xi\deriv{}{x}$ est
nilpotente donc pas semi-simple.

On note $\A$ la sous-algèbre de Lie de $\Cinf(M)$
fonctionnellement engendrée par les $f_1,\dots,f_n$. Si
 $g_1,\dots,g_n$ est un système générateur de $\A$, de telle
sorte qu'il existe des fonctions $F_1,\dots,F_n$ telles que 
\[ f_i = F_i(g_1,\dots,g_n),~ i=1,\dots,n \]
alors les $\He(g_i)$ forment une base de $A$ qui s'obtient à partir de
 $\He(f_i)$ par le changement de base linéaire $dF^{-1}(0)$.
Ceci montre en particulier que la sous-algèbre de Cartan associée ne
dépend pas de la base de $\A$ choisie.

On a alors le théorème fondamental~:
\begin{theo}[\mbox{Eliasson \cite[4.1]{eliasson-these}}]
\label{TEliasson} 
Pour un tel système complètement in\-tég\-rab\-le, soit $(q_1,\dots,q_n)$
une base standard de l'algèbre de Cartan associée. Alors le feuilletage
lagrangien singulier donné par les surfaces de niveau des $f_i$
est localement symplectiquement égal à celui donné par les $q_i$.

Autrement dit, il existe un \diffeo\  symplectique $\phy$ au voisinage
de $m$ tel que~:
\[ \forall i,j~~ \{ f_i\circ\phy,q_j\}=0. \]
\end{theo}

\subsection{Un exemple~: le problème de C.Neumann classique}
\label{neumann}

Le problème de C.Neumann est celui du mouvement d'une particule sur
une sphère de dimension $n$ soumise à une force dérivant d'un potentiel
quadratique.  On se donne donc une matrice symétrique réelle $A$
définie positive de taille $n+1$. Le potentiel $V$ est la
restriction à $S^n$ de la forme quadratique sur $\RM^{n+1}$ associée~:
\[ V(x) = \frac{1}{2}\pscal{Ax}{x}. \]
Le mouvement est décrit sur $T^*S^n$ par l'Hamiltonien 
\[H(x,\xi) = \frac{1}{2}|\xi|^2 + V(x).\] La métrique sur la sphère est celle
induite par la métrique euclidienne sur $\RM^{n+1}$, et la structure
symplectique sur $T^*S^n$ n'est autre que la restriction de la
structure standard de $T^*\RM^{n+1}$.  Remarquons ici que le problème
est invariant par antipodie, ce qui nous autorise à le
considérer sur $\PM^n$ plutôt que sur $S^n$.

Il est pratique de voir ce système comme la restriction d'un système
hamiltonien sur $T^*\RM^{n+1}$ de la façon suivante~:
 
 Contraindre le mouvement à s'effectuer sur $S^n$ revient à tenir
compte d'une force de réaction normale à la sphère. L'équation du
mouvement prend la forme~:
\[  \ddot{x}=-V'(x) + \lambda x. \]
De $|x|^2=1$, on tire
$\pscal{x}{\dot{x}}=\pscal{x}{\ddot{x}}+|\dot{x}|^2=0$, ce qui permet
de déterminer $\lambda$~:
\[ \lambda = \pscal{V'(x)}{x} - |\dot{x}|^2 = 2V(x)- |\dot{x}|^2.\]
On peut alors vérifier que cette équation est obtenue en restreignant
à \begin{equation}
\label{cotangent}
 T^*S^n = \{(x,\xi), ~~ |x|^2=1, ~~ \pscal{x}{\xi} = 0\}
  \end{equation}
 le champ de vecteurs sur $T^*\RM^{n+1}$ d'Hamiltonien~:
\[ H_0(x,\xi) = V(x) + \frac{1}{2}(|x|^2|\xi|^2 - \pscal{x}{\xi}^2). \]

On peut alors montrer que $H_0$ est complètement intégrable, et on
dispose même d'un système explicite d'intégrales en involution, trouvé
en 1975 par Uhlenbeck. Ce problème a des liens très étroits avec le
flot géodésique sur un $n$-ellipsoïde; à ce sujet, voir
\cite{moser-ens}.
 
Nous nous intéressons ici aux points fixes du flot.
\begin{prop}
  Le problème de C.Neumann sur $T^*\PM^n$ a exactement $n+1$ points
  fixes $p_0,\dots,p_n$. En ordonnant convenablement les $p_i$, pour
  un potentiel non-résonnant, le système complètement intégrable associé
  est, au voisinage de $p_i$, non-dégénéré de type $(n-i,i,0)$.
\end{prop}
\demo on détermine les points critiques de $H$ par la méthode des
multiplicateurs de Lagrange, qui se trouvent être les valeurs propres
$\lambda_i$ de $A$. Les $p_i$ sont donc une base de vecteurs propres
aux voisinages desquels on peut écrire la forme de Morse de $H$. On
voit alors que si $V$ est assez générique, $H''$ est un élément
régulier de la sous-algèbre de Cartan, ce qui permet de déterminer son
type. \cqfd

\rem[1] la condition de non-résonance requise est celle de
l'indépendance sur $\ZM$, pour tous $i$, des
$\sqrt{|\lambda_i-\lambda_j|}$, $j\neq i$. Elle est bien-entendu
générique.

\rem[2] on n'obtient pas dans cet exemple de singularité
\emph{focus-focus}. Mais ce dernier type de singularité, même s'il a
longtemps été négligé, est loin d'être rare; il apparaît par exemple
au point d'équilibre instable du pendule sphérique (voir
\cite{duistermaat}).

\subsection{Le commutant classique}
De la même façon qu'au niveau quadratique l'objet essentiel était pour
nous la sous-algèbre de Cartan $\mathfrak{c}_F$, on voit à l'énoncé du
théorème \ref{TEliasson} que l'algèbre de fonctions $\Cinf$ qui va
nous intéresser tout au long de ce travail est le \emph{commutant
classique} des $q_j$, qu'on note $C_q$~:
\[ C_q = \{f\in \Cinf(\RM^{2n}),  \forall j,~ \{f,q_j\}=0 \}. \]
C'est l'algèbre des fonctions localement constantes sous l'action des
champs Hamiltoniens des $q_i$. 

Décrivons la structure de cette algèbre. Une première idée naturelle
est que, au moins formellement, tout élément de $C_q$ doit pouvoir
s'exprimer comme ``fonction des seules variables $q_1,\dots,q_n$''. En
dehors du cadre formel, ce n'est cependant pas aussi simple.

Notons $z_i=(x_i,\xi_i)$, et soit $H_i$ l'ensemble des points
$(z_1,\dots,z_n)\in\RM^{2n}$ tels que $z_i=0$. L'union des $H_i$ est
exactement le lieu des points critiques de l'application moment
$q=(q_1,\dots,q_n)$.
\begin{lemm}
  Soit $U$ un ouvert de $\RM^{2n}\setminus (\bigcup_{i=1}^n H_i)$ tel
  que pour tous $c\in\RM^n$, la feuille lagrangienne $U\cap q^{-1}(c)$
  soit connexe (ou vide). Alors la restriction de toute fonction $f\in
  C_q$ à $U$ peut s'écrire comme fonction $\Cinf$ de $(q_1,\dots,q_n)$
  uniquement.
\end{lemm}
\demo en chaque point $z$ de $\RM^{2n}\setminus (\bigcup_{i=1}^n
H_i)$, les fonctions $q_1,\dots,q_n$ sont indépendantes, donc, d'après
le théorème de Darboux-Carathéodory, elles peuvent être complétées en
des coordonnées symplectiques $(y_1,\dots,y_n,q_1,\dots,q_n)$ sur un
voisinage $\Omega$ de $z$. La condition d'appartenance à $C_q$ devient
alors $\deriv{f}{y}=0$. En recouvrant $U$ par de tels voisinages
$\Omega_\alpha$, on obtient des fonctions $F_\alpha$ telles que
\[ f_{\restr \Omega_\alpha} = F_\alpha(q_1,\dots,q_n). \] L'hypothèse
de connexité des surfaces de niveau de $q$ implique que la valeur de
$F_\alpha(q_1,\dots,q_n)$ ne dépend pas de $\alpha$. On obtient ainsi
une unique fonction $\Cinf$ $F$ telle que 
\[ f_{\restr U} = F(q_1,\dots,q_n), \]
ce qui répond à la question. \cqfd

Une conséquence directe de lemme est que $C_q$ est une
algèbre \emph{commutative}. En effet, deux fonctions
$F_1(q_1,\dots,q_n)$ et $F_2(q_1,\dots,q_n)$ commutent. Donc deux
éléments $f_1$ et $f_2$ de $C_q$ commutent presque partout,
et donc partout, par continuité du crochet de Poisson.

Supposons maintement que $n=1$ et que l'algèbre de Cartan
$\mathfrak{c}_q$ soit de type $(1,0,0)$ (autrement dit, $q$ est
elliptique). Les surfaces de niveau de $q$ sont des cercles; on peut
donc, dans le précédent lemme, choisir $U=\RM^2\setminus\{0\}$ et
obtenir que toute $f\in C_q$ s'écrive $F(q)$ sur $U$. On peut
montrer facilement qu'un tel $F$ est en réalité $\Cinf$ en $0$, et
donc que $f$ s'écrit $F(q)$ sur tout $\RM^2$. Un même phénomène se
produit lorsque $n=2$ et que $\mathfrak{c}_q$ est de type focus-focus
$(0,0,1)$; il suffit pour cela de faire le calcul en coordonnées
complexes
\[ u=x_1+\sqrt{-1}x_2,  \quad v=\xi_1+ \sqrt{-1}\xi_2,\]
 dans lesquelles les surfaces de niveau s'écrivent $\bar{u}v=cste$ et
sont connexes (voir \cite{eliasson-these} pour le détail des calculs).

Le cas pathologique est donc celui d'une algèbre de type $(0,1,0)$
(hyperbolique). En effet les surfaces de niveau ont chacune deux
composantes connexes. Un élément de $C_q$ est donc décrit
par deux fonctions $F_+$ et $F_-$ telles que, par exemple,
\begin{equation}
  \label{recolle}
  f_{\restr x>0} = F_+(q), \quad f_{\restr x<0} = F_-(q).
\end{equation}
On voit facilement que $F_+-F_-$ doit être plate en $0$, et
réciproquement, tout couple de fonctions $(F_+,F_-)$ sur $\RM$ tel que
$F_+-F_-$ est plat en $0$ définit un élement de $C_q$ par
(\ref{recolle}).

On obtient ainsi l'énoncé suivant~:
\begin{prop}
\label{commutant-classique}
$C_q$ est une algèbre de Poisson commutative. Si la
sous-algèbre de Cartan $\mathfrak{c}_q$ est de type $(m_e,m_h,m_l)$,
on note $(q_1,\dots,q_n)$ une base standard ordonnée, comportant en
premier lieu les éléments hyperboliques.  Pour chaque $m_h$-uplet
$\epsilon=(\epsilon_1,\dots,\epsilon_{m_h})$, $\epsilon_i=\pm 1$, on
pose \[ E_\epsilon=\{(x,\xi)\in\RM^{2n},~\forall i=1,\dots,m_h, \quad
\epsilon_i x_i >0 \}. \] Un élément $f$ de $C_q$ est
caractérisé par la donnée de $2^{m_h}$ fonctions
$F_\epsilon\in\Cinf(\RM^n)$ telles que~:
\begin{itemize}
\item $\forall \epsilon$, $f_{\restr E_\epsilon} =
  F_\epsilon(q_1,\dots,q_n)$.
\item $\forall (\epsilon,\epsilon')$, $F_\epsilon-F_{\epsilon'}$ est
  plate sur $\overline{E_\epsilon}\cap \overline{E_{\epsilon'}}$.
\end{itemize}
\end{prop}
On en déduit un ``lemme de division'' que nous
utiliserons souvent~:
\begin{lemm}
  \label{taylor}
  Soit $g\in C_q$. Alors 
  \begin{itemize}
  \item il existe des fonctions $g_i\in C_q$ telles que~:
    \[ g = g(0) + \sum_{i=1}^n g_i q_i. \]
  \item il existe des   fonctions $a_i\in C_q$ telles que~:
    \[ \ham{g} = \sum_{i=1}^n a_i \ham{q_i}. \]
  \end{itemize}
\end{lemm}
\demo par la formule de Taylor, chaque $G_\epsilon$ correspondant à
$g$ et donné par la proposition \ref{commutant-classique} s'écrit:~
\[ G_\epsilon(q_1,\dots,q_n) = G_\epsilon(0) +  \sum_{i=1}^n q_i
G_\epsilon^i(q_1,\dots,q_n). \] Pour chaque $i$, les couples
$(G_\epsilon^i,G_{\epsilon'}^i)$ vérifient les mêmes hypothèses de
platitude, et donc définissent un élément $g_i$ de
$C_q$. Ces $g_i$ répondent au premier point.

De la même façon, sur chaque $E_\epsilon$ on a
\[ \ham{g} = \sum_{i=1}^n \deriv{G_\epsilon}{q_i} \ham{q_i}. \]
En posant $a_\epsilon^i=\deriv{G_\epsilon}{q_i}$, on voit que les
$a_\epsilon^i$ vérifient toujours les bonnes hypothèses de
platitude. Ils définissent donc des fonctions $a_i\in C_q$
qui répondent au deuxième point. \cqfd

Revenant au théorème d'Eliasson, la proposition
\ref{commutant-classique} implique immédiatement l'existence de formes
normales ``fortes'' en l'absence d'élément hyperbolique~:
\begin{coro}
\label{TFN}
Soit $\A$ un système complètement intégrable critique non dé\-gé\-né\-ré;
soit $(q_1,\dots,q_n)$ une base standard de l'algèbre de Cartan
associée. On suppose qu'il n'y a pas parmi les $q_i$ d'élément
hyperbolique. Alors il existe des coordonnées symplectiques dans
lesquelles $\A$ est engendrée par les $q_i$.

En d'autres termes, il existe un \diffeo\  symplectique $\phy$ au
voisinage de $m$ et des fonctions $F_i~:~\RM^n\fleche\RM$, de
différentielles indépendantes en 0, telles que~:
\[ \forall i=1,\dots,n~~ f_i\circ\phy = F_i(q_i,\dots,q_n). \]
\end{coro}

\rem[1] il n'est pas certain que l'hypothèse de connexité des feuilles
lagrangiennes soit nécessaire, puisque le résultat a été prouvé pour
toute singularité non dégénérée en dimension deux~: c'est le ``Lemme
de Morse isochore'' (voir \cite{colin-vey}).  En outre, on sait
(``théorème de Sternberg'') qu'en dimension quelconque $n$, un
hamiltonien $H$ hyperbolique non résonnant (c'est-à-dire dont la
hessienne est régulière dans une sous-algèbre de Cartan totalement
hyperbolique: elle s'écrit $\sum \lambda_i x_i\xi_i$ avec des
$\lambda_i$ indépendants sur $\ZM$) est complètement intégrable~: il
existe un symplectomorphisme local $\phy$ tel que $H\circ\phy$ soit
une fonction des $x_i\xi_i$. Or si l'on se donne un système
complètement intégrable $\mathcal{A}$ de type purement hyperbolique
($m_h=n$), il existe une base de $\mathcal{A}$ formée par des
fonctions toutes non-résonnantes. Il ne paraît donc pas absurde de
penser qu'il puisse exister une forme normale dans ce cas. Dans
\cite{guillemin-schaeffer} est exposée une belle démonstration du
théorème de Sternberg mais elle semble résister au cas complètement
intégrable...

Cependant, bien que le résultat serait intéressant en lui-même, une
forme normale générale ne nous serait pas utile puisque le lemme
crucial \ref{poincare}, nécessaire pour traiter le problème au niveau
semi-classique est, quant à lui, faux en présence d'un mélange
d'éléments hyperboliques avec des éléments d'autres types.

\rem[2] lorsque la singularité est entièrement \emph{elliptique}
 ($q_i=x_i^2+\xi_i^2$), le corolaire \ref{TFN} est montré d'une façon
 différente par \cite{dufour-molino}.

\rem[3] dans la catégorie analytique, la structure du commutant perd
toute sa richesse, puisqu'il n'y a pas de fonctions analytiques
plates, en dehors des fonctions constantes... D'ailleurs l'analogue
analytique du théorème \ref{TEliasson} était connu bien avant (voir
\cite{russmann} en dimension 2, et \cite{vey} en dimension
quelconque).

Si maintenant $\mathfrak{c}_q$ contient des éléments hyperboliques
($m_h\neq 0$), nous utiliserons l'énoncé suivant, corollaire du
théorème \ref{TEliasson} et du lemme \ref{taylor}~:
\begin{coro}
\label{TFN-faible}
Soient $f_1,\dots,f_n$ des fonctions définissant un système
comp\-lè\-te\-ment in\-tég\-rab\-le critique non dégénéré, et
$q_1,\dots,q_n$ une base standard de l'algèbre de Cartan associée (de
type quelconque); alors il existe un symplectomorphisme $\phy$ tel
qu'on ait, au voisinage de $0$,~:
\[ (f_1,\dots,f_n)\circ \phy = M.(q_1,\dots,q_n), \]
où $M$ est une matrice $n\times n$ de fonctions éléments de
$C_q$, inversible en $z=0$.
\end{coro}
Le fait que $M$ soit inversible en $0$ est exactement l'hypothèse de
non-dé\-gé\-né\-res\-cence du système~: $M(0)$ est la matrice telle
que~:
\[ \left( \He(f_1),\dots,\He(f_n)\right) =
M(0).\left(\He(q_1),\dots,\He(q_n)\right) \]
Notons que, sur l'ouvert où $M$ est inversible, son inverse reste
formé d'élé\-ments de $C_q$, puisque ce dernier est une
algèbre de Poisson.

\vspace{0.5cm}

\subsection{Un ``Lemme de Poincaré'' critique}
\label{Spoincare}

Après les théorèmes de forme normale \ref{TEliasson} et \ref{TFN},
l'ingrédient le plus important dont nous auront besoin pour mener à
bien la quantification semi-classique est le résultat suivant~:

\begin{theo}[Eliasson]
\label{poincare}
Soit $q_1,\dots,q_n$ une base standard d'une sous-algèbre de Cartan
$\mathfrak{c}_q$ de $\Lie{sp}(2n)$. Dans tout voisinage de $0$, il existe un
sous-voisinage $\Omega$ de $0$, tel que, si $g_1,\dots,g_n$ sont des
fonctions $\Cinf$ à valeurs réelles ou complexes vérifiant~:
\begin{equation}
  \label{hypothese}
  \forall i,j=1,\dots,n \quad \{g_i,q_j\} = \{g_j,q_i\} \textrm{ sur } \Omega, 
\end{equation}
alors il existe une fonction $f$ définie sur $\Omega$, et des fonctions
$F_i\in C_q$  formellement uniques, telles que, sur $\Omega$,
\begin{equation}
\label{systeme}
 \forall i=1,\dots,n ~~  \{f,q_i\} = g_i - F_i.
\end{equation}
\end{theo}
En réalité, l´énoncé  d'Eliasson dans \cite{eliasson-these} ne
mentionne pas l'existence d'un voisinage universel $\Omega$, mais
celui-ci se déduit facilement de sa preuve. On rappellera comment lors
de la preuve de la variante \ref{poincare-hyperbolique}.

\rem[1] pour $c\in\RM^n$, notons $\Omega_c$ la lagrangienne $\cap
q_i^{-1}(c_i)$. Elle est singulière pour $c=0$. Aux points réguliers,
les champs hamiltoniens $\ham{q_i}$ forment une base de l'espace
tangent à $\Omega_c$, ce qui permet de définir une 1-forme $g_c$ sur
$\Omega_c$ par~:
\[ g_c (\ham{q_i}) = g_i. \]
L'hypothèse du théorème \ref{poincare} traduit alors que cette 1-forme
est fermée. L'in\-té\-grer serait trouver une fonction $f$ vérifiant
$ \{f,q_i\} = g_i $. Le théorème montre qu'on peut trouver une
primitive $f$ dépendant régulièrement du paramètre $c$ même en $c=0$,
pourvu qu'on retranche à $g$ une certaine fonction de $c$ uniquement.
Cette fonction contient donc les singularités ``gênantes'' de $g$ en
$0$.

\rem[2] la remarque précédente prouve que l'hypothèse \ref{hypothese}
du théorème est nécessaire, ce qui se vérifie aussi immédiatement en
utilisant l'identité de Jacobi.

Supposons maintenant que $f_1,\dots,f_n$ soient des fonctions
définissant un système comp\-lè\-te\-ment in\-tég\-rab\-le critique
non dégénéré en un point $m$, et soit $(q_1,\dots,q_n)$ une base
standard du la sous-algèbre de Cartan associée. Appliquant le théorème
\ref{TEliasson}, on peut toujours supposer que, localement, les $f_i$
sont dans $C_q$. Le théorème \ref{poincare} admet alors le
corollaire suivant~:
\begin{coro}
\label{corollaire-poincare-faible}
Avec les hypothèse précédentes, dans tout voisinage de $m$, il existe
un sous-voisinage $\Omega$ tel que, si $g_1,\dots,g_n$ sont des
fonctions vérifiant~: \[ \forall i,j=1,\dots,n ~~ \{g_i,f_j\} =
\{g_j,f_i\} \textrm{ sur } \Omega, \] alors il existe une fonction $a$
sur $\Omega$ et des fonctions $F_i$ dans $C_q$ telles que
\begin{equation}
 \forall i=1,\dots,n ~~  \{a,f_i\} = g_i - F_i.
\end{equation}
\end{coro}
\demo
une application du lemme \ref{taylor} donne l'existence d'une matrice
$N$ de fonctions dans $C_q$, inversible en $0$, telle que~:
\[ (\ham{f_1},\dots,\ham{f_n}) = N.(\ham{q_1},\dots,\ham{q_n}). \]
La 1-forme $g_c$ définie sur les parties régulières des feuilles
lagrangiennes $\Omega_c$ par $g_c(\ham{f_i})=g_i$, est donc donnée,
dans la base $(\ham{q_1},\dots,\ham{q_n})$, par
\[ (g_c(\ham{q_1}),\dots,g_c(\ham{q_n}) = N^{-1}.(g_1,\dots,g_n) \]
ce qui implique que l'hypothèse de fermeture $\{g_i,f_j\} =
\{g_j,f_i\}$ est équivalente à $\{\tilde{g}_i,q_j\} =
\{\tilde{g}_j,q_i\}$, où on a noté 
\[(\tilde{g}_1,\dots,\tilde{g}_n)=N^{-1}.(g_1,\dots,g_n). \]

Appliquant alors le théorème \ref{poincare}, on obtient, sur $\Omega$,
une fonction $a$ et des fonctions $\tilde{F}_i$ telles que~:
\[ (\{a,q_1\},\dots,\{a,q_n\}) =  (\tilde{g}_1,\dots,\tilde{g}_n) -
(\tilde{F}_1,\dots,\tilde{F}_n), \] ce qui se réécrit~:
\[ N.(\{a,q_1\},\dots,\{a,q_n\}) = (g_1,\dots,g_n) -
N.(\tilde{F}_1,\dots,\tilde{F}_n) .\] Or
$N.(\{a,q_1\},\dots,\{a,q_n\}) = (\{a,f_1\},\dots,\{a,f_n\})$, donc
les fonctions $a$ et 
\[ (F_1,\dots,F_n) =
N.(\tilde{F}_1,\dots,\tilde{F}_n)\] répondent à la question. \cqfd

Bien-sûr, diverses variantes de l'énoncé du théorème \ref{poincare} et
du corollaire ci-dessus s'obtiennent immédiatement en utilisant la
proposition \ref{commutant-classique} et le lemme \ref{taylor}.  En
particulier, lorsque la sous-algèbre de Cartan $\mathfrak{c}_q$ ne
contient pas d'élément hypoerbolique ($m_h=0$), les $F_i$ s'écrivent
comme fonctions $\Cinf$ de $(q_1,\dots,q_n)$

Nous allons montrer maintenant que tel est aussi le cas dans la
situation extrême où $\mathfrak{c}_q$ est de type entièrement hyperbolique
($m_h=n$).
\begin{theo}
  \label{poincare-hyperbolique}
  Avec les hypothèses du théorème \ref{poincare}, si en outre la
  sous-algèbre de Cartan $\mathfrak{c}_q$ est de type hyperbolique
  $(0,n,0)$, alors il existe une fonction $f$ définie sur $\Omega$, et
  des fonctions $\Cinf$ $F_1,\dots,F_n$ de $\RM^n$ dans $\RM$ ou $\CM$
  telles que, sur $\Omega$,
  \begin{equation}
    \label{systeme-hyperbolique}
    \forall i=1,\dots,n ~~  \{f,q_i\} = g_i - F_i(q_1,\dots,q_n).
  \end{equation}
\end{theo}
On obtient ainsi un résultat optimal, dans le sens où une telle
formulation devient fausse dès $m_h\neq 0$ et $m_h\neq n$.

\demo tout au long des sections (\ref{cas-hyp-debut}) et
(\ref{cas-hyp-fin}), on pose donc $q_i=x_i\xi_i$.  Le schéma de la
preuve est classique~: on commence par examiner le cas formel
(\ref{cas-hyp-debut}), qui résout le problème modulo des fonctions
plates, puis la solution du système (\ref{systeme-hyperbolique}) est
donnée par une formule explicite (\ref{cas-hyp-fin}).

\subsubsection{Le cas hyperbolique formel}
\label{cas-hyp-debut}
On note $\mathcal{P}_k(\RM^{2n})$ l'espace vectoriel des polynômes homogènes
de degré $k$ en les variables $x_1,\dots,x_n,\xi_1,\dots,\xi_n$.
L'algèbre de Lie (pour le crochet de Poisson) des fonctions formelles
est constituée par l'espace
\[ \mathcal{P} = \bigoplus_{k=0}^\infty \mathcal{P}_k. \]
Si $P\in\mathcal{P}_m$, la re\-pré\-sen\-ta\-tion adjointe
\[ \ad{P}~:~Q\mapsto \{P,Q\} \]
est un endomorphisme de $\mathcal{P}$ envoyant $\mathcal{P}_k$ sur
$\mathcal{P}_{m+k-2}$. En particulier, si $P$ est quadratique,
$\ad{P}$ est un endomorphisme de chaque $\mathcal{P}_k$, et 
$\ad{}$ est une re\-pré\-sen\-ta\-tion de ${\mathcal Q}(2n)=\mathcal{P}_2$.

Soit $\hat{\A}$ la sous-algèbre des fonctions formelles des $q_i$
(une base de $\hat{\A}$ est ainsi constituée par les monômes 
$x^{\al}\xi^{\beta}$
avec $\al = \beta$). Par abus de notation, on identifiera parfois
$\hat{\A}$ avec $\mathcal{P}(\RM^{n})$ par $F=F(q_1,\dots,q_n)$.

On veut donc prouver le résultat suivant~:
\begin{prop}
  Soient $g_1,\dots,g_n \in \mathcal{P}$ telles que
  \[ \forall i,j=1,\dots,n ~~ \{g_i,q_j\} = \{g_j,q_i\}. \]
  Alors il existe une fonction $f\in\mathcal{P}$, et des fonctions
  $F_i\in\hat{\A}$ telles que
  \[ \forall i=1,\dots,n ~~  \{f,q_i\} = g_i - F_i(q_1,\dots,q_n). \]
\end{prop}
\demo
pour tous $i$, $\hat{\A}\in\ker(\ad{q_i})$ donc $\ad{q_i}$ se quotiente en
un endomorphisme de $\mathcal{P}/\hat{\A}$
 (on désigne ainsi le quotient gradué 
$\oplus\mathcal{P}_k/\hat{\A}_k$). 
Plaçons-nous pour toute cette
démonstration dans $\mathcal{P}/\hat{\A}$; la proposition s'énonce alors
ainsi~:
\begin{equation}
\label{sysformel}
\exists f,~~ \forall i,~~ \ad{q_i} f = g_i.
\end{equation}
Une algèbre de Cartan peut être caractérisée par un élément régulier,
aussi n'est-il pas étonnant que le système (\ref{sysformel}) puisse,
comme on va le voir, se réduire à une seule équation, à savoir
\[ \ad{q} f = g, \textrm{ où } q=\sum \la_i q_i \textrm{ et } g=\sum \la_i g_i, \]
où les $\la_i$ sont fixés indépendants sur $\ZM$.

 $\ad{q}$ étant diagonalisable, de valeurs propres
 $\sum \la_i(\beta_i -\al_i)$ (associées aux vecteurs propres 
 $x^{\al}\xi^{\beta}$), elle est inversible dans 
$\mathrm{End}(\mathcal{P}/\hat{\A})$ 
( = $\oplus\mathrm{End}(\mathcal{P}_k/\hat{\A}_k)$ )
et $f$ est donc uniquement déterminée (dans $\mathcal{P}/\hat{\A}$).

Montrons que $f=(\ad{q})^{-1}g$ convient~:
\[ \ad{q_i} f = \ad{q_i} (\ad{q})^{-1}g =  (\ad{q})^{-1} \ad{q_i} g \]
car $\ad{q_i}$ et $ \ad{q}$ commutent (car $\{q_i,q\} = 0$).
Or l'hypothèse $\{g_i,q_j\} = \{g_j,q_i\}$ entraîne 
\[  \ad{q_i} g =  \ad{q} g_i \]
donc $ \ad{q_i} f = g_i$ et la proposition est démontrée. \cqfd

\subsubsection{Preuve du théorème \ref{poincare-hyperbolique}}
\label{cas-hyp-fin}
En appliquant la proposition précédente aux jets d'ordre infini
$\hat{g_i}$ de $g_i$ (qui vérifient toujours $\{\hat{g_i},q_j\} =
\{\hat{g_j},q_i\}$), et en prenant des représentants $\Cinf$ des
fonctions $\hat{f}$ et $F_i$ obtenues, on obtient que la fonction
\[  \{\hat{f},q_i\} - ( \hat {g_i} - F_i(q_1,\dots,q_n)) \]
est plate en $0$.
Posons $f=\hat{f}+h$; on est ainsi ramené à chercher une fonction $h$
telle que $\{h,q_i\} = \tilde{g_i}$, où $\tilde{g_i}=g_i - F_i - 
\{\hat{f},q_i\}$ est plate en $0$. Remarquons que $\tilde{g_i}$
vérifie toujours l'hypothèse \ref{hypothese}.
%

Il suffit donc de montrer~:
\begin{prop}
\label{Pplat}
Dans tout voisinage de $0$, il existe un voisinage $\Omega$ de $0$ tel
que, si $g_i$ sont des fonctions $\Cinf$ plates en $0$, vérifiant, sur
$\Omega$,
 \[ \{g_i,q_j\} = \{g_j,q_i\}, \]
alors il existe une fonction $f$ sur $\Omega$ (plate en $0$) telle que
 \[ \forall i=1,\dots,n ~~  \{f,q_i\} = g_i. \]
\end{prop}
\demo il serait tentant d'utiliser les méthodes de
\cite{guillemin-schaeffer}.  Malheureusement, cela nécessiterait de
pouvoir supposer que les $g_i$ soient à support compact, ce qui ne
paraît pas compatible avec (\ref{hypothese}). Toujours est-il qu'en
dimension 1, l'hypothèse (\ref{hypothese}) est vide, et la méthode
marche très bien (voir \cite{guillemin-schaeffer} (theo 4)).

Nous allons donc nous rabattre sur les techniques d'intégration standard,
expliquées en dimension 1 dans \cite{colin-vey} et
\cite{eliasson-these}. On montrera
ici comment elles permettent de traiter le cas complètement intégrable.

\begin{figure}[h]
    \begin{center}
    \leavevmode
    \input{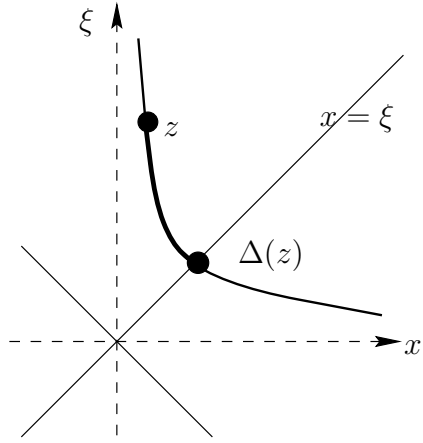}
    \caption{$\Delta(z)=U(T(z))z$.}
    \label{fig:T}
    \end{center}
  \end{figure}

On désigne toujours par $U_i(t)$ les flots des $q_i$. Notons $T_i(z)$
le temps que met $z$ sous l'action de $U_i$ à placer la coordonnée
$z_i$ sur l'unique hyperplan diagonal ou antidiagonal $x_i=\pm\xi_i$
qui rencontre le flot partant de $z$. C'est
une fonction bien définie et $\Cinf$ en dehors des hyperplans de
coordonnées. On montre dans les références sus-citées que, pour toute
fonction $g$ plate en $0$, l'intégrale~:
\[ \int_0^{T_i(z)}g(U_i(s)z)ds \]
définit une fonction $\Cinf$ sur un voisinage de $0$.  Notons
$\Delta_i(z)$ la projection de $z$ sur le $i$ème hyperplan
(anti)diagonal le long du flot de $q_i$, c'est-à-dire
$\Delta_i(z)=U_i(T_i(z))z$ (voir figure \ref{fig:T}). Il est essentiel
de remarquer ici que tout voisinage de $0$ contient un sous-voisinage
stable par $\Delta_i$; on choisit pour $\Omega$ un tel voisinage stable.
\begin{lemm}
\label{enfin}
Soit $g$ une fonction plate en $0$ et, pour $z\in\Omega$, soit
$f_i(z)= \int_0^{T_i(z)}g(U_i(s)z)ds$. Alors~:
\[ \{f_i,q_i\} = g, \] et, pour $j\neq i$, si $h$ est une fonction
vérifiant $\{g,q_j\} = \{h,q_i\}$ sur $\Omega$, alors~:
\[ \{f_i,q_j\} = h - \Delta_i^*h. \] 
\end{lemm}
\demo c'est un simple calcul, reposant sur~:
\[ \{f_i,q_j\} = -{\ddt}_{t=0} f_i(U_j(t)). \]
Pour $i=j$, on utilise que $T_i(U_i(t)z)=T_i(z)-t$, tandis que pour
$j\neq i$, c'est l'invariance de $T_i$ par $U_j$, jointe à l'hypothèse
de commutation, qui donne le résultat. Cette hypothèse est utilisée
sous la forme:~
\[ \ddt g(U_j(t)U_i(s)z) = \{q_j,g\}U_j(t)U_i(s)z = \]
\[ = \{q_i,h\}U_i(s)U_j(t)z = \frac{d}{ds} h(U_i(s)U_j(t)z). \]
~ \cqfd

Revenons à la preuve de la proposition. En appliquant le lemme avec
$i=j=1$, on résout la première équation: $\{f_i,q_i\} = g_1$.  On
cherche alors une solution du système de la forme $f_1+f_2$, ce qui
mène aux équations~:
\[ \{f_2,q_1\} = 0, \textrm{ et }\forall j>1, ~~ \{f_2,q_j\} =
\tilde{g_j} \] où $\tilde{g_j} = g_j - \{f_1,q_j\}$. Le membre de
droite vaut $\Delta_1^*g_j$, d'après le lemme. Les fonctions
$\Delta_1^*g_j$ et $T_j$ sont invariantes par le flot de $q_1$, donc
une nouvelle application du lemme avec $g=\Delta_1^*g_2$ résout la
deuxième équation tout en laissant la première intacte. On peut ainsi
recommencer... jusqu'à épuisement du système.  Le lemme assurera que
les dérivées croisées sont les bonnes.

On peut même s'offrir une formule explicite, dont on vérifie la
validité grâce au lemme~:
\[  f = \int_0^{T_1(z)} g_1(U_1(s)z) ds +
\int_0^{T_2(z)} \Delta_1^*g_2(U_2(s)z) ds + \]
\[ +
\int_0^{T_3(z)} \Delta_2^*\Delta_1^*g_3(U_3(s)z) ds + \cdots \]
~ \cqfd

\section{Systèmes intégrables quantiques}
\label{quantique}

\subsection{Introduction}
D'une certaine façon, on peut dire que l'intégrabilité d'un point de
vue purement quantique ne présente guère d'intérêt. Pour voir cela,
prenons $(E,h\omega)$ un espace vectoriel symplectique (de
dimension finie, au moins dans un premier temps). La matrice
symplectique standard $J$, vérifiant $J^2=-I$, munit $E$ d'une
structure complexe.  D'autre part $E$ est muni d'un produit scalaire
euclidien par $(x,y)=\omega(Jx,y)$, qui n'est autre que le produit
euclidien usuel.  On obtient ainsi une structure hermitienne sur $E$
grâce au produit hermitien
 $\pscal{x}{y} = (x,y) +ih\omega(x,y)$. 
 Réciproquement, tout espace de Hilbert
 $\CM^n$ est muni de la forme symplectique égale à la partie imaginaire
de la forme hermitienne.

On vérifie alors facilement qu'un endomorphisme de $E$ est
\emph{hermitien} si et seulement s'il commute avec $J$ et est
symétrique (pour le produit scalaire euclidien). Ou plutôt, 
voyant $\mathrm{Herm}(E,J)$ comme $J\Lie{u}(n)$, on a~:
\[ \mathrm{Herm}(E,J) = J(\Lie{sp}(E) \cap \Lie{o}(E)). \]

On veut maintenant résoudre l'équation de Schrödinger associée à $B$,
à savoir~:
\[ \frac{h}{i}\ddt \psi = B\psi, \]
où on a noté $i$ pour $J$. Si $H(\psi) = (B\psi,\psi)$ est la forme
quadratique associée à $B$, on vérifie immédiatement que
l'hamiltonien de $H$ est $\ham{H}(\psi) = \frac{i}{h}B(\psi)$. L'équation de
Schrödinger n'est donc autre que l'équation du flot de   $\ham{H}$~:
\[\ddt \psi =\ham{H}(\psi).\]
Or un hamiltonien hermitien est évidemment complètement intégrable,
puisqu'en le diagonalisant en base hermitienne, on obtient une base
symplectique (et orthonormée) dans laquelle
\[ B=\left(\begin{array}{cc}\Delta & 0 \\ 0 &
\Delta\end{array}\right), ~~~ \Delta = 
\mathrm{diag}(\lambda_1,\dots,\lambda_n). \]
 $B$ s'écrit donc $\sum \lambda_i(x_i^2+\xi_i^2)$ et les
 $x_i^2+\xi_i^2$ sont des intégrales en involution.

On peut remarquer ici que le crochet de Poisson de deux hamiltoniens
hermitiens s'écrit 
\[ \{B_1,B_2\} = \frac{i}{h}[B_1,B_2], \]
où $[,]$ est le crochet usuel des endomorphismes.

Maintenant, sans rentrer dans des justifications mathématiques, si
$\mathcal{H}$ est un espace de Hilbert quelconque sur lequel on se
donne un opérateur $B$ autoadjoint diagonalisable, avec des fonctions
propres $\psi_n$ associées aux valeurs propres $\lambda_n$, alors
l'équation de Schrödinger peut être vue comme un système hamiltonien
de dimension infinie, d'Hamiltonien $H(\psi) = \sum \lambda_n
|z_n|^2$, ($z_n=x_n + i\xi_n$ est la coordonnées de $\psi$ sur
$\psi_n$) pour la structure symplectique $h\sum d\xi_n \wedge dx_n$.
Il est donc ``complètement intégrable'', au sens où il admet les
$x_n^2 + \xi_n^2$ comme famille complète d'intégrales premières.

Ceci dit, on n'est guère plus avancé... C'est pourquoi on introduit
une autre notion de complète intégrabilité quantique, qu'on pourra
aussi appeler la complète intégrabilité semi-classique, puisqu'elle va
coïncider avec la notion classique dans la limite semi-classique.

\subsection{Quantification semi-classique}
On note $\mathcal{E}$ l'espace des $h$-\opds\ classiques d'ordre 0 sur
un ouvert $W$ d'une variété $X$ de dimension $n$. Rappelons qu'un
\opd\ $P(h)$ est dit \emph{classique} si son symbole de Weyl $p(h)$
admet un développement asymptotique de la forme~:
\[ p(h)(x,\xi) \sim \sum h^k p_k(x,\xi) \]
(le développement est alors unique). Les résultats présentés ici
 seront microlocaux au voisinage d'un point ou d'une sous-variété
 Lagrangienne compacte; en particulier, peu importent les hypothèses
 que l'on peut faire sur le comportement à l'infini des symboles (voir
 \cite{courscolin,robert}). Rappelons aussi que deux \opds\  sont dits
 microlocalement égaux (ou égaux modulo $O(h^\infty)$) au voisinage
 d'un point $m\in T^*X$ si leurs symboles de Weyl coïncident sur ce
 voisinage.

$\mathcal{E}$ est, comme en dimension finie, muni d'une structure
 d'algèbre de Lie~:
\[ \{P,Q\} = \frac{i}{h}[P,Q]. \] L'opération $\sigma$ (symbole
principal) est un morphisme surjectif de $(\mathcal{E},\{\})$ sur
l'algè\-bre des observables classiques $(\Cinf(T^*W),\{\})$. 
\begin{defi}
Suivant \cite{colinII,charbonnel}, on appelle \emph{système
  semi-classique complè\-te\-ment intég\-rab\-le} la donnée de $n$
  opérateurs $P_1,\dots,P_n$ de $\mathcal{E}$ tels que~:
  \begin{itemize}
  \item les symboles principaux $p_i=\sigma(P_i)$ sont des fonctions à
    valeurs réelles formant un système complètement intégrable
    classique;
  \item $\forall i,j$, $[P_i,P_j]=O(h^\infty)$. 
  \end{itemize}
\end{defi}
Cette notion est bien adaptée à ce qu'il semble naturel d'appeler une
``quantification semi-classique'' d'un système classique complètement
intégrable.

Par exemple, le problème de Neumann présenté en section \ref{neumann}
est semi-classiquement intégrable; il est même \emph{quantiquement}
intégrable, au sens où il existe des opérateurs différentiels sur
$S^n$ qui quantifient les Hamiltoniens classiques en commutant
exactement (voir \cite{gurarie}).

Cette définition soulève naturellement la question suivante~:


Lorsqu'un système complètement intégrable classique associé à un
système complètement intégrable semi-classique est réductible à une
forme normale (par exemple, s'il satisfait aux hypothèses du théorème
\ref{TFN}; c'est aussi vrai dès que l'application moment
$J=(p_1,\dots,p_n)$ est une submersion locale), existe-t'il aussi une
forme normale pour le système semi-classique ?

La réponse à cette question est connue, et positive, au voisinage de
points réguliers des feuilles lagrangiennes du système. Elle est
principalement due à Colin de Verdière \cite{colinII} et
A.-M. Charbonnel \cite{charbonnel} (même si, rigoureusement, ces
auteurs ne traitent que le cas où les opérateurs sont formellement
auto-adjoints).

Nous verrons dans la section suivante que cela reste le cas en
présence d'une singularité non-dégénérée~: tous les résultats énoncés
en section \ref{classique} auront leur analogue semi-classique.

%
%
%
%

\subsection{Formes normales microlocales}

Soit $(P_1(h),\dots,P_n(h))$  un système semi-classique complètement
intégrable, défini sur un ouvert $W$ d'une variété $X$ de dimension
$n$. Rappelons que les $P_i$ ne sont pas nécessairement
formellement auto-adjoints, mais que le symbole principal est supposé
réel.

On note $p_j$ les symboles principaux; ils définissent un système
complètement intégrable classique de $T^*W$. On suppose ici que ce
système admet un point critique non-dégénéré $m\in T^*W$, et on note
$\mathfrak{c}_q$ la sous-algèbre de Cartan de $\Lie{sp}(2n)$ associée.
On suppose toujours que $p_j(m)=0$.

Soit $(q_1,\dots,q_n)$ une base standard de $\mathfrak{c}$. On note
$Q_j$ les $h$-quantifiés sy\-mét\-ri\-ques (ou de Weyl) des $q_j$,
c'est-à-dire~:
\begin{itemize}
        \item si $q_j = x_j^2+\xi_j^2$, ~~
        $Q_j = -h^2 \deriv{\,^2}{x_j^2} + x_j^2$;
        \item si $q_j = x_j\xi_j$, ~~
        $Q_j = \frac{h}{i}(x_j\deriv{}{x_j} + \frac{1}{2})$;
        \item et pour le cas focus-focus, on a respectivement
        $Q_j=\frac{h}{i}(x_j\deriv{}{x_{j+1}} -x_{j+1}\deriv{}{x_j})$ et
        $Q_{j+1}=\frac{h}{i}(1 +x_j\deriv{}{x_{j}} +x_{j+1}\deriv{}{x_{j+1}})$
(ce qui donne en coordonnées polaires $\left\{\begin{array}{ll}
 x_j=\rho\cos(\theta) \\ x_{j+1}=\rho\sin(\theta) \end{array}\right.$, 
les opérateurs respectifs $\frac{h}{i}\deriv{}{\theta}$ et 
 $\frac{h}{i}(\rho\deriv{}{\rho}+1)$).
\end{itemize}



Outre la sous-algèbre de Cartan $\mathfrak{c}_q$, l'ingrédient
principal de la section \ref{classique} était une sous-algèbre de
Poisson de $\Cinf(T^*W)$, le commutant classique $C_q$. On définit de
la même façon le \emph{commutant semi-classique} $\mathcal{C}_Q$ comme
étant la sous-algèbre de $\mathcal{E}$ donnée par~:
\[ \mathcal{C}_Q = \{ P(h)\in \mathcal{E}, \forall j,~
[P,Q_j]=O(h^\infty)\}. \] 

La version semi-classique du théorème d'Eliasson est alors la
suivante~:
\begin{theo}
  \label{TQ3}
  Soit $(P_1,\dots,P_n)$ un système semi-classique comp\-lè\-te\-ment
  in\-té\-gra\-ble, avec une singularité non dégénérée de type
  quelconque, et $Q_1,\dots,Q_n$ les opé\-ra\-teurs différentiels
  correspondants.  Alors il existe un \opif\ elliptique (et unitaire
  si les $P_j$ sont formellement auto-adjoints) $U(h)$, et des \opds\
  $F_j(h)$ dans $\mathcal{C}_Q$, tels que, microlocalement au voisinage
  de $m$,
\[ \forall j ,~~ U^{-1}P_jU = F_j + O(h^\infty). \]
\end{theo} 

Avant de donner la preuve de ce théorème, notre première tâche est
donc de décrire la structure de $\mathcal{C}_Q$.

\subsection{Le commutant semi-classique}
$\mathcal{C}_Q$ est directement relié au commutant classique $C_q$ par
le résultat suivant~:
\begin{prop}
\label{commutant-quantique}
$\mathcal{C}_Q$ est une algèbre de Lie commutative. L'application
 $\sigma_W$ (symbole de Weyl) est, modulo $O(h^\infty)$, un
 isomorphisme de $\mathcal{C}_Q$ dans l'algèbre de Lie formelle~:
\[ C_q(h) = \sum_{k=0}^\infty h^k C_q. \]
\end{prop}
\demo la formule de Moyal (voir \cite{flato}) exprime la relation
entre la quantification de Weyl et le crochet de commutation~: si $A$
et $B$ sont les quantifiés de Weyl des symboles $a$ et $b$, on a,
formellement~:
\[ \sigma_{Weyl}( [A,B]) = \frac{2}{i}a\sin(\frac{h\mathcal{D}}{2}) b,
\]
avec \[ \mathcal{D} = \left(\overleftarrow{\deriv{}{\xi}}
\overrightarrow{\deriv{}{x}} - \overleftarrow{\deriv{}{x}}
\overrightarrow{\deriv{}{\xi}}\right). \]
On l'utilise en géneral pour la formule suivante~:
\[ \sigma_{Weyl}( [A,B]) = \frac{h}{i}\left( \{a,b\} + O(h^2)\right).
\]

Commençons par montrer la commutativité de $\mathcal{C}_Q$.  Sur
chaque ouvert $E_\epsilon$ (voir proposition
\ref{commutant-classique}), les symboles principaux des $Q_j$ sont
indépendants. D'après le théorème de forme normale standard associé à
une carte de Darboux-Ca\-ra\-théo\-do\-ry $(x,\xi)$ (voir par exemple
\cite{colinII}), les $Q_j$ sont conjugués aux opérateurs
$\frac{h}{i}\deriv{}{x_j}$. En appliquant la formule de Moyal, on
obtient que tout \opd\ qui commute avec les
$\frac{h}{i}\deriv{}{x_j}=Op^W(\xi_j)$ a un symbole de Weyl dont le
dé\-ve\-lop\-pe\-ment asymptotique est indépendant de $x$. Par une
nouvelle application de cette formule on voit que deux tels opérateurs
commutent à l'ordre $O(h^\infty)$.

Finalement, si $P(h)$ et $Q(h)$ sont dans $\mathcal{C}_Q$, le crochet
$[P,Q]$ est un \opd\ classique dont le symbole de Weyl a un
développement asymptotique nul sur les ouverts $E_\epsilon$. Par
continuité des termes de ce développement asymptotique, ils sont nuls
partout. Donc $[P,Q]=O(h^\infty)$.

Venons-en maintenant au deuxième point de la proposition.  Appliquant
la formule de Moyal dans les coordonnées initiales, le fait que $q_j$
soit quadratique implique que pour tous $k\geq 3$ et pour toutes
fonctions $f$,
\[ f\mathcal{D}^k q_j =  0. \]
On dispose donc d'une règle de quantification ``exacte'', au sens où~:
\[\sigma_{Weyl}( [Op^W(f),Q_j]) =  \frac{h}{i}\left( \{f,q_j\} +
  O(h^\infty) \right),\] On en déduit qu'un opérateur $P(h)$ est dans
$\mathcal{C}_Q$ si et seulement si le dé\-ve\-lop\-pe\-ment
asymptotique de son symbole de Weyl $p(h)\sim\sum p_kh^k$ vérifie
\[ \sum h^k\{p_k,q_j\} = O(h^\infty), \]
ce qui est équivalent, par unicité du développement asymptotique, à~:
\[ \forall k,j, \quad  \{p_k,q_j\}=0,\]
c'est-à-dire $\forall k$, $p_k\in C_q$. \cqfd

Lorsque la sous-algèbre de Cartan $\mathfrak{c}_q$ ne contient pas
d'élément hyperbolique ($m_h=0$), on sait que la structure de $C_q$
est particulièrement simple~: c'est l'ensemble des fonctions
$f(q_1,\dots,q_n)$. On a un résultat analogue pour le commutant
semi-classique.

Notons tout d'abord que puisque les opérateurs $Q_j$ sont formellement
auto-adjoints, on peut utiliser un calcul fonctionnel à plusieurs
variables comme celui développé dans \cite{charbonnel}. On peut ainsi
donner un sens microlocal à $f(Q_1,\dots,Q_n)$ pour toute fonction
$f\in\Cinf_0(\RM^n)$ et par là, à $f(h)(Q_1,\dots,Q_n)$, où $f(h)\sim
f_0+hf_1 + ...$ est un symbole semi-classique admettant un
développement en puissance de $h$~: $f(h)(Q_1,\dots,Q_n)$ commute alors
avec les $Q_j$ et admet $f_0(q_1,\dots,q_n)$ comme symbole principal.

Il faut remarquer que microlocalement sur un domaine $\Omega$
voisinage de $m$, les hypothèses à l'infini pour le calcul fonctionnel
sont inutiles. Ainsi, seule la symétrie des opérateurs est
importante. En effet, quitte à tronquer les symboles en dehors de
$\Omega$, on peut supposer que les opérateurs $h$-quantifiés sont
semi-bornés inférieurement. On dispose alors par exemple de
l'extension de Friedrichs qui est auto-adjointe.  Enfin, deux
extensions d'un même opérateur, dont les domaines contiennent
$\Cinf_0(\Omega)$, auront le même symbole de Weyl sur $\Omega$, et
seront donc égaux modulo régularisant. En particulier, le choix d'une
extension auto-adjointe est sans incidence sur le calcul fonctionnel,
au niveau microlocal.

On peut maintenant énoncer la proposition suivante~:
\begin{prop}
  \label{fonctionnel}
  Si $m_h=0$, alors tout élément de $\mathcal{C}_Q$ s'écrit,
  microlocalement au voisinage de $0$, sous la forme
  $f(h)(Q_1,\dots,Q_n) + O(h^\infty)$.
\end{prop}
\demo soit $P(h)\in\mathcal{C}_Q$. D'après la proposition
\ref{commutant-quantique}, son symbole principal appartient à $C_q$,
et donc s'écrit $f^{(0)}(q_1,\dots,q_n)$. Posons, microlocalement au
voisinage de $0$,
\[ P^{(0)}=f^{(0)}(Q_1,\dots,Q_n).\]
$P$ et $P^{(0)}$ ayant même symbole principal, on a~:
\[ P = P^{(0)} + hR^{(1)}, \]
où, nécessairement, $R^{(1)} \in \mathcal{C}_Q$. En recommençant la
même décomposition pour $R^{(1)}$, et ainsi de suite, on obtient par
récurrence~:
\[ \forall N, \exists f^{(0)}, \dots,f^{(N)} \in \Cinf(\RM^n), \quad
P(h) = (\sum_{k=0}^N h^kf^{(k)})(Q_1,\dots,Q_n) + O(h^{N+1}), \]
ce qui prouve la proposition. \cqfd

\subsection{Preuve du théorème \ref{TQ3}}
Plutôt que de montrer directement le théorème \ref{TQ3}, nous allons
énoncer et prouver une formulation légèrement différente qui se
rapproche du corollaire \ref{TFN-faible}, et qui est plus utile pour
les applications lorsque $m_h\neq 0$. Dans le cas $m_h=0$, la forme la
plus utile est bien-sûr celle du théorème \ref{TQ3} associé à la
proposition \ref{fonctionnel}.
\begin{theo}
  \label{TQ2} 
  Soit $(P_1,\dots,P_n)$ un système semi-classique comp\-lè\-te\-ment
  in\-té\-grable, avec une singularité non dégénérée de type quelconque,
  et $Q_1,\dots,Q_n$ les opérateurs différentiels correspondants.
  Alors il existe un \opif\ elliptique (et unitaire si les $P_j$ sont
  formellement auto-adjoints) $U(h)$, une matrice de taille $n\times
  n$ microlocalement inversible $\mathcal{M}(h)$ d'\opds\ appartenant
  à $\mathcal{C}_Q$, et des constantes $\alpha_j^{(\ell)}\in\CM$
  ($j=0,\dots,n$ et $\ell\in\NM^*$) telles que, microlocalement au
  voisinage de $m$,
  \[ U^{-1}(P_1,\dots,P_n)U =
  \mathcal{M}.(Q_1-\alpha_1(h),\dots,Q_n-\alpha_n(h)) +
  O(h^\infty). \] On a noté $\alpha_j(h)=\sum_{\ell\geq
  1}\alpha_j^{(\ell)}h^\ell$.

  Les $\alpha_j^{(1)}$ se décrivent à l'aide des valeurs en $m$ des
  symboles sous-principaux $r_j$ des $P_j$ par~:
  \[ \alpha^{(1)}=-M^{-1}(0).r(m),\] où $M$ est la matrice de fonctions
  donnée par le corollaire \ref{TFN-faible}.
\end{theo}
\demo On commence par appliquer le corollaire \ref{TFN-faible}. En
choisissant un \opif\ microlocalement unitaire $U$ associé au
symplectomorphisme $\phy$, et en conjuguant les $P_j$ par $U$, on se
ramène au cas où les symboles principaux $p_j$ vérifient, dans un
voisinage de $0$,
\[ (p_1,\dots,p_n)=M.(q_1,\dots,q_n) . \] 
$M$ étant une matrice de fonctions commutant avec les $q_k$, une
application de la proposition \ref{commutant-quantique} nous fournit
une matrice $\mathcal{M}^{(0)}$ d'\opds\  éléments de $\mathcal{C}_Q$,
de symbole principal $M$. En conséquence, il existe des opérateurs
$R_j^{(1)}\in\mathcal{E}$ tels que~:
\[ (P_1,\dots,P_n) =  \mathcal{M}^{(0)}.(Q_1,\dots,Q_n) +
h(R^{(1)}_1,\dots,R^{(1)}_n). \] Notons
$(M_1,\dots,M_n)=\mathcal{M}^{(0)}.(Q_1,\dots,Q_n)$. Les $M_j$ sont
des éléments de $\mathcal{C}_Q$, de symboles principaux
$(p_1,\dots,p_n)$.  L'hypothèse de commutation des $P_j$ s'écrit~:
\[ [M_j,M_k] + h\left([M_j,R_k]+[R_j,M_k]\right) + h^2[R_j,R_k] = 0. \] 
Puisque $\mathcal{C}_Q$ est abélien (proposition
\ref{commutant-quantique}), on obtient, en prenant les symboles
principaux de cette égalité~:
\[ \forall j,k, \quad \{p_j,r_k\} = \{p_k,r_j\}.\]
La preuve du théorème se poursuit alors en deux étapes~: la première
traitant spécifiquement du niveau \emph{sous-principal}, et la
deuxième traitant, par un même schéma, tous les niveaux suivants.

$\bullet$ La première étape consiste à chercher une matrice
$\mathcal{M}^{(1)}$ d'éléments de $\mathcal{C}_Q$, des constantes
$\alpha_j^{(1)}\in\CM$, et un \opd\ elliptique $V\in\mathcal{E}$ tels
que~:
\[ V^{-1}(P_1,\dots,P_n)V =
(\mathcal{M}^{(0)} + h\mathcal{M}^{(1)}) .
(Q_1-h\alpha_1^{(1)},\dots,Q_n-h\alpha_n^{(1)})+ O(h^2). \] 
Ce système se réécrit en~:
\[ \left[\sum_{k=1}^n\mathcal{M}^{(0)}_{jk}Q_k,V\right] = h\left(V\sum_{k=1}^n
  \mathcal{M}^{(1)}_{jk}Q_k -
  V\sum_{k=1}^n\mathcal{M}^{(0)}_{jk}\alpha_k^{(1)} - R_jV\right) +
O(h^2), \]
qui est équivalent, en prenant les symboles principaux des deux
  membres, au système d'équations de transport suivant~:
\[ \forall j, \quad \{p_j,v\} = iv\left(\sum_{k=1}^n m^{(1)}_{jk}q_k -
  \sum_{k=1}^n m^{(0)}_{jk}\alpha_k^{(1)} - r_j\right).\]
Posons $v=e^{id}$; on a  $\{f,v\} = iv\{f,d\}$; il suffit donc de résoudre~:
\[ \forall j, \quad \{p_j,d\} = \left(\sum_{k=1}^n m^{(1)}_{jk}q_k -
  \sum_{k=1}^n m^{(0)}_{jk}\alpha_k^{(1)} - r_j\right).\] Le
corollaire \ref{corollaire-poincare-faible} fournit une telle fonction
$d$ définie sur un voisinage $\Omega$ de $0$ et des éléments $F_j$ de
$C_q$ tels que
\[  \forall j, \quad \{p_j,d\} = F_j-r_j. \]
La matrice $(m^{(0)}_{jk})=M$ étant inversible, le n-uplet
$M^{-1}.(F_1,\dots,F_n)$ est le plus général des n-uplets de $C_q^n$,
donc, d'après le lemme \ref{taylor}, il s'écrit
\[ -(\alpha_1^{(1)},\dots,\alpha_n^{(1)}) + \tilde{M}.(q_1,\dots,q_n)
,\]
où $\tilde{M}$ est une matrice d'éléments de $C_q$, et les
$\alpha_j^{(1)}$ sont nécessairement donnés par~:
\[ \forall j,~~ -\sum_{k=1}^n m^{(0)}_{jk}(0)\alpha_k^{(1)} = F_j(0) =
r_j(0).\]
Le système est donc résolu, avec $(m^{(1)}_{jk})=M.\tilde{M}$.

En quantifiant les $m_{jk}^{(1)}$ par la proposition
\ref{commutant-quantique}, on obtient une matrice $\mathcal{M}^{(1)}$
d'\opds\  de $\mathcal{C}_Q$ qui résout le problème modulo $O(h^2)$.

Notons que si les $P_j$ sont formellement auto-adjoints, alors $d$ est
réel, et $V$ peut être quantifié en un opérateur microlocalement
unitaire.

$\bullet$ La deuxième étape termine la preuve par récurrence~:
supposons que les $P_j$ vérifient
\[ (P_1,\dots,P_n) =  \left(\sum_{\ell=0}^{N-1}
  h^\ell\mathcal{M}^{(\ell)}\right).(Q_1-\alpha_1(h),\dots,Q_n-\alpha_n(h))
+ h^{N}(R^{(N)}_1,\dots,R^{(N)}_n), \] où 
$ \alpha_j(h) = \sum_{\ell=1}^{N-1} h^\ell\alpha_j^{(\ell)}$,
et les $\mathcal{M}^{(\ell)}$ sont des matrices $n\times n$ d'éléments
de $\mathcal{C}_Q$. Comme tous ces éléments commutent entre eux,
l'hypothèse de commutation des $P_j$ s'écrit, à l'ordre principal
$N$~:
\[ h^N \left( [M_j,R^{(N)}_k] + [R^{(N)}_j,M_k] \right) =
  O(h^\infty), \]
ce qui donne~:
\[ \{p_j,r_k^{(N)}\} = \{p_k,r_j^{(N)}\}. \]
Alors il existe un \opd\ $C\in\mathcal{E}$, une matrice
$\mathcal{M}^{(N)}$ d'élé\-ments de $\mathcal{C}_Q$, et des constantes
$\alpha_j^{(N)}\in\CM$ telles que~:
\begin{eqnarray}
  \label{etape2}
  \lefteqn{\forall j, \quad
    (I+h^{N-1}iC)^{-1}(P_1,\dots,P_n)(I+h^{N-1}iC) = {} } \nonumber \\
    & & \left(\sum_{\ell=0}^{N}
    h^\ell \mathcal{M}^{(\ell)}\right).(Q_1- \tilde{\alpha}_1(h),
    \dots,Q_n- \tilde{\alpha}_n(h)) + O(h^{N+1}),
\end{eqnarray}
avec $\tilde{\alpha}_j(h)= \alpha_j(h) + h^N\alpha_j^{(N)}$.  En
effet, pour $N\geq 2$, ce système se réécrit modulo des termes d'ordre
$O(h^{N+1})$ en~:
\[ h^{N-1}\left[\sum_{k=1}^n\mathcal{M}^{(0)}_{jk}Q_k,iC\right] =
h^N\left(\sum_{k=1}^n \mathcal{M}^{(N)}_{jk}Q_k -
  \sum_{k=1}^n\mathcal{M}^{(0)}_{jk}\alpha_k^{(N)} -
  R^{(N)}_j\right). \] Il donne lieu, au niveau des symboles
principaux, à des équations de transport que l'on résout comme
précédemment.

Notons que si les $P_j$ sont formellement auto-adjoints, le symbole
$c$ est réel. Au lieu de conjuguer par
l'opérateur elliptique $I+h^{N-1}iC$, on utilise
la transformée de Cayley
$W=\frac{I+ih^{N-1}C/2}{I-ih^{N-1}C/2}$. L'opérateur $W$ est unitaire,
et comme $W=I+ih^{N-1}C+O(h^N)$, il convient tout aussi bien pour la
résolution de (\ref{etape2}).

Ceci termine la démonstration du théorème \ref{TQ2}. \cqfd

Le théorème \ref{TQ3} en est bien-sûr un corollaire immédiat.

\subsection{Remarques finales}
\paragraph{1.}
 Le schéma de démonstration du théorème \ref{TQ2} est général, et
fonctionne si l'on remplace l'algèbre $\mathcal{C}_Q$ par celle des
fonctions des $(Q_1,\dots,Q_n)$, à partir du moment où les théorèmes
classiques le permettent. Ceci est légitime si la sous-algèbre de
Cartan ne contient pas d'élément hyperbolique -- on l'a d'ailleurs
déjà remarqué; ça l'est aussi en dimension 1, en vertu du lemme de
Morse isochore (\cite{colin-vey}), et du lemme de Poincaré
hyperbolique \ref{poincare-hyperbolique}.  Cette remarque permet de
retrouver le théorème 12 de \cite{colin-p}.

\paragraph{2.}
 On peut aussi utiliser ce schéma de démonstration pour avoir
la forme normale des $(P_1,\dots,P_n)$ dans un voisinage invariant
d'une feuille \emph{régulière} du feuilletage lagrangien déterminé par
l'application moment $(p_1,\dots,p_n)$. En effet le théo\-rème
d'Arnold-Liouville fournit un tel voisinage invariant $\Omega$ dans
lequel on dispose de coordonnées action-angle, c'est-à-dire que
$\Omega$ est symplectomorphe à un voisinage de la section nulle de
$T^*\T^n$ de la forme $\T^n\times\mathcal{D}\ni(x,\xi)$. Utilisant un
\opif\  associé à ce symplectomorphisme, on est ramené au cas où les
$P_j$ agissent sur $\T^n$ et ont des symboles principaux qui ne
dépendent que de $\xi$.

Examinons maintenant les équations de transport~: elles vont consister
à ré\-sou\-dre des systèmes du type 
\[ \{ \xi_j,a\} = r_j, \textrm{ c'est-à-dire }\deriv{a}{x_j} = r_j. \]
Or les $r_j=r_j(x,\xi)$ vérifient la bonne condition de fermeture~:
\[  \{\xi_j,r_k\} = \{\xi_k,r_j\},\]
qui n'est autre que $d_x r = 0$, ou $r$ est vue comme la 1-forme 
$r_1dx_1+\cdots +r_ndx_n$. Cette condition assure qu'on peut toujours
intégrer ces équations localement. Par contre, pour un résultat global
sur les tores horizontaux, il faut $[r]=0$ dans $H^1(\T^n\times\{\xi\})$.
Autrement dit, on sait résoudre $d_x a = r - [r]$, ou encore~:
\[ \{ \xi_j,a\} = r_j-[r]_j, \]
où $[r_j]$ est la moyenne de $r_j$ par rapport à $x_j$, et ne dépend
que de $\xi$.

On retrouve ainsi le résultat mentionné en introduction sur les
coordonnées actions-angles semi-classiques régulières~:
\begin{theo}[\cite{colinII,charbonnel}]
  \label{TQreg}
  Soit $(P_1,\dots,P_n)$ un système semi-classiquement
  comp\-lè\-te\-ment intégrable, et $\Omega$ un ouvert invariant de
  coordonnées actions-angles ré\-gu\-li\-ères. Alors il existe un
  \opif\ $U(h)$ (unitaire si les $P_j$ sont formellement
  auto-adjoints) et des symboles
  \[ f_j(h) = f_j^{(0)} + h f_j^{(1)} + h^2 f_j^{(2)} + \cdots \]
  tels que, sur $\Omega$,
  \[ \forall j = 1,\dots,n, ~~~ U^{-1}P_jU = f_j(h)(D_1,\dots,D_n) +
  O(h^\infty). \]
  Les $D_j$ sont les opérateurs $\frac{h}{i}\deriv{}{x_j}$ sur le
  tore $\T^n$. 
\end{theo}

\paragraph{3.}
 Dans l'énoncé du théorème \ref{TQ2}, la constante $\alpha^{(1)}$ est
définie de façon intrinsèque par les opérateurs $P_j$, c'est-à-dire
indépendemment de l'OIF $U$ choisi. Il semblerait raisonnable de
penser qu'il en est de même de tous les termes de la série. Un
argument en faveur de cette affirmation est que $\alpha(h)$ peut être
vue comme la monodromie des solutions du système pseudo-différentiel
$P_ju=0$ (voir \cite{pham,san2}). Nous espérons pouvoir détailler cela
dans un travail ultérieur.

\noindent \textbf{Remerciements}\\
Je tiens à remercier Yves Colin de Verdière pour ses remarques
pertinentes.  

\newpage

\bibliographystyle{plain}
\bibliography{bibli}


\end{document}